\documentclass[12pt]{amsart}
\usepackage{color,amssymb,amsmath,amsthm,amsfonts,mathrsfs,enumerate,esint,MnSymbol,latexsym}
\usepackage{a4wide}
\usepackage[english]{babel}
\theoremstyle{plain}

\numberwithin{equation}{section}

  \newcommand{\const}{\rm const}

  \newcommand{\grad}{\rm grad}

 \DeclareMathOperator*{\esssup}{ess\,sup}

\newtheorem{theorem}{Theorem}[section]
\newtheorem{definition}{Definition}[section]

\newtheorem{remark}{Remark}[section]

\renewenvironment{proof}{{\bf{Proof.}}}{\hfill $\Box$ \\}

\newenvironment{nouppercase}{%
  \renewcommand{\uppercasenonmath}[1]{}}{}

\begin{document}

\title[FUNDAMENTAL SOLUTION FOR CAUCHY INITIAL VALUE PROBLEM FOR
PARABOLIC PDEs]{Fundamental solution for Cauchy initial value
problem for parabolic PDEs with discontinuous unbounded first-order
coefficient at the origin.\\ Extension of the classical parametrix
method}

\footnotesize\date{}

\author[Maria Rosaria Formica, Eugeny Ostrovsky and Leonid Sirota]{Maria Rosaria Formica ${}^{1}$, Eugeny Ostrovsky ${}^2$, Leonid Sirota ${}^2$}

\begin{nouppercase}
\maketitle
\end{nouppercase}

\begin{center}
${}^{1}$ Parthenope University of Naples, via Generale Parisi 13,\\
Palazzo Pacanowsky, 80132,
Napoli, Italy. \\

e-mail: mara.formica@uniparthenope.it \\

\vspace{4mm}

${}^2$ Department of Mathematics and Statistics, Bar-Ilan University, \\
59200, Ramat Gan, Israel. \\

e-mail: eugostrovsky@list.ru\\

e-mail: sirota3@bezeqint.net\\

\end{center}

\begin{abstract}

We prove the existence of a fundamental solution of the Cauchy
initial boundary value problem on the whole space for a parabolic
partial differential equation with discontinuous unbounded
first-order coefficient at the origin. We establish also
non-asymptotic, rapidly decreasing at infinity, upper and lower
estimates for the fundamental solution. We extend the classical
parametrix method provided by E.E. Levi.
\end{abstract}

%

\vspace{2mm}

 \noindent {\footnotesize {\it Key words and phrases}:

Partial Differential Equation of parabolic type, parametrix method,
fundamental solution, degenerate diffusion, Gamma and Beta-
functions, Markov random process, transfer density of probability
for diffusion random processes, infinitesimal operator, generalized
Mittag-Leffler function, Chapman-Kolmogorov equation, Neumann
series, H\"older's continuity, Volterra's integral equation.

\vspace{2mm}


\vspace{5mm}

\section{Definitions. Notations. Statement of problem.}

\vspace{5mm}

We consider in this article the Partial Differential Equation (PDE)
of (uniform) parabolic type

\begin{equation} \label{equation1}
\frac{\partial u}{\partial t} = \frac{1}{2} \frac{\partial^2 u}{\partial x^2} + \frac{b(t,x)}{|x|^{\gamma}} \ \frac{\partial u}{\partial x},
\end{equation}
where $u = u(t,x)$, $ t \ge 0, \ x \in \mathbb R, \  \ \gamma =
\const \in (0,1)$ and  the coefficient $b = b(t,x)$ is a continuous
bounded function such that $ \ b(t,0) \ne 0, \ $  together with the
initial value problem (Cauchy statement)

\begin{equation} \label{initcond}
\lim_{t \to s+} u(t,x) = f(x), \ \ u = u(t,x) = u[f](t,x),  \ \ s
\ge 0,
\end{equation}
 where $f(\cdot)$ is a measurable function satisfying some growth condition at  $ |x| \to \infty$;  for the definiteness one can suppose its continuity and
 boundedness, i.e.
$\esssup_x |f(x)| < \infty$. \par
 \ Briefly, $ \ L_{t,x}u = L[b]_{t,x}u= 0, \ $ where $ \ L[b]_{t,x}  \ $ is the linear parabolic partial differential operator of the form
\begin{equation} \label{operator L}
L[b]_{t,x} = L_{t,x} \stackrel{def}{=} \frac{\partial }{\partial t} - \frac{1}{2} \frac{\partial^2 }{\partial x^2} - \frac{b(t,x)}{|x|^{\gamma}} \ \frac{\partial} {\partial x}.
\end{equation}

 \vspace{4mm}

Let us recall the classical definition.

\vspace{4mm}

\begin{definition}
{\rm The (measurable) function $ \ p = p(t,x,s,y), \ 0 < s < t; \
x,y \in \mathbb R^2 $ is said to be a \emph{Fundamental Solution}
(F.S.) for the equation \eqref{equation1}, subject to the initial
condition \eqref{initcond}, iff
 for all the fixed values $ \ (s,y), \ 0 < s < t, \ x,y \in  \mathbb R^2 \ $ it satisfies the equation \eqref{equation1} and, for all the bounded continuous functions $f =
 f(x)$,
\begin{equation} \label{funsol}
\lim_{t \to s+} \int_{\mathbb R} p(t,x, s,y) \ f(y) \ dy  = f(x), \
\ x \in  \mathbb R,
\end{equation}
if, of course, there exists.}
\end{definition}

In this case the F.S. $p(t,x,s,y)$ may be interpreted as the
transfer density of probability for diffusion random non-homogeneous
in the time, as well as in the space, Markov process, see e.g.
\cite{Bally}, \cite{Corielli}.  As a consequence:

$$
p(t,x, s,y)  \ge 0; \ \ \int_{\mathbb R} p(t,x, s,y) \ dy = 1
$$
and satisfies the well-known Chapman-Kolmogorov equation

\begin{equation} \label{Chapman Kolm}
p(t,x,s,y) = \int_{\mathbb R} p(t,x,v,z) \ p(v,z,s,y) \ dz, \ 0 < s
< v < t.
\end{equation}

 \vspace{4mm}

 On the other words, the (unique) solution $u = u(t,x)$ of the Cauchy problem \eqref{equation1}, \eqref{initcond}, in the whole space $\mathbb
 R$,
 may be written as follows

 \begin{equation} \label{solution}
 u(t,x) = \int_{\mathbb R} p(t,x,s,y) \ f(y) \, dy.
 \end{equation}

  \vspace{3mm}

 The solution of the non-homogeneous parabolic PDE of the form

\begin{equation} \label{equation non homogen}
\frac{\partial u}{\partial t} = \frac{1}{2} \frac{\partial^2
u}{\partial x^2} + \frac{b(t,x)}{|x|^{\gamma}} \ \frac{\partial
u}{\partial x} + g(t,x), \ \ t > s,
\end{equation}
with zero initial condition $u(s,x) = 0$, may be written, as
ordinary, as follows

\begin{equation} \label{solutionnonhom}
 u(t,x) = \int_s^t d \theta \int_{\mathbb R} p(t,x,\theta,y) \ g(\theta,y) \ dy.
 \end{equation}

\vspace{4mm}

 For instance, the well-known Gauss function

\begin{equation} \label{Zfun}
Z(t,x,s,y) = Z_{t-s}(x-y) := (2 \pi (t-s))^{-1/2} \ \exp \left\{ \ -
\frac{1}{2} \ \frac{(x-y)^2}{(t-s)}  \ \right\},
\end{equation}

$0 < s < t < \infty, \ x,y \in \mathbb R$, is the Fundamental
Solution (F.S.) for the ordinary heat equation

$$
\frac{\partial u}{\partial t}  =  \frac{1}{2} \ \frac{\partial^2 u}{\partial x^2}.
$$

 \vspace{5mm}

Our claim is to prove the existence of F.S. for the equation
\eqref{equation1} with the initial condition \eqref{initcond}, under
suitable conditions.

\vspace{5mm}

 The existence of F.S. for {\it strictly}  parabolic PDE with smooth coefficients, for instance, bounded and belonging to the H\"older space, is  explained in many works, e.g.
 \cite{Bally}, \cite{Dressel2}, \cite{DeckKruse}, \cite{Friedman}, \cite{Gagliardo}, \cite{IlKalashOleynik}, \cite{Ladyzhenskaya}, \cite{Lax}, \cite{Levi1}, \cite{Levi2}.
The same problem for parabolic systems is considered in
\cite{Eydelman1}, \cite{Eydelman2}, \cite{Slobodetskiy}. The case of
parabolic equations with discontinuous but bounded coefficients is
investigated in \cite{DeckKruse}, \cite{Ilyin}, \cite{Kamynin},
\cite{Krasnosel}, \cite{Samarskiy}, \cite{Sobolevskiy} etc.\par
 We will  essentially follow the interesting article of Thomas Deck and Susanne Kruse \cite{DeckKruse}, where the authors considered the case of the equation

\begin{equation} \label{equation3}
\frac{\partial u}{\partial t} = \frac{1}{2} \frac{\partial^2 u}{\partial x^2} + b_1(t,x) \ \frac{\partial u}{\partial x},
\end{equation}
in which the coefficient $ \ b_1 = b_1(t,x) \ $ is H\"older
continuous in $x$ and may be unbounded at $ \ x \to \pm \infty$, in
the following sense
$$
|b_1(t,x) | \le C \  \left( \ 1 + |x|^{\beta} \ \right), \ \
0<C,\beta = \const < \infty.
$$

 The uniqueness of the F.S. for parabolic PDEs is proved, e.g. in \cite{Dressel1},  \cite{Dressel2}, \cite{Gagliardo}, \cite{Kamynin}, \cite{Ladyzhenskaya}, \cite{Samarskiy}. \par

\vspace{4mm}

In the probability theory, more precisely, in the theory of
diffusion random processes, the F.S. $ p(t,x,s,y) $ represents the
density of transition of probability. The equations of the form
\eqref{equation1} appears in particular in \cite{Molchanov},
\cite{Paul Eric}; see also \cite{Bally}, \cite{Corielli}.

\vspace{5mm}

\section{Main result: Existence of the Fundamental solution. Upper and lower estimates.}

\vspace{5mm}

 \begin{theorem}\label{Th FS}
Suppose that the coefficient $ b = b(t,x)$ in the equation
\eqref{equation1} is measurable and
 bounded, i.e.
\begin{equation} \label{bound}
 K = K[f] := \sup_{t,x} |b(t,x)| = ||b(\cdot,\cdot)||_{C((0,\infty) \otimes \mathbb R)} <
 \infty.
\end{equation}
Then there exists a Fundamental Solution $ p(t,x,s,y) $ for the
equation \eqref{equation1}.

 \end{theorem}

\vspace{0.3cm}

\begin{proof}

\vspace{0.3cm}

 \ {\bf A.}  We will follow the proof in \cite{DeckKruse}, herewith we retain
analougus notations, so that we state that the F.S. is given by
\begin{equation} \label{Urav pt}
p(t,x,s,y) = Z_{t-s}(x-y) + \int_s^t \int_R Z_{t-r}(x-z) \  \Phi(r,z,s,y) \ dr \ dz,
\end{equation}
where the function  $ \ \Phi = \Phi(t,x,s,y) \ $  is a solution of the following Volterra equation

\begin{equation} \label{Volterr Phi}
\Phi(t,x,s,y) = \Phi_1(t,x,s,y) + \int_s^t \int_R \Phi_1(t,x,r,z) \ \Phi(r,z,s,y) \ dr \ dz
\end{equation}
with weak singular  kernel

\begin{equation} \label{Expr for Phi1}
\Phi_1(t,x,s,y) = \frac{b(t,x)}{|x|^{\gamma}} \ \frac{\partial }{\partial x} \ Z_{t-s}(x-y).
\end{equation}

 \vspace{3mm}
The solution of the equation \eqref{Volterr Phi} may be obtained, as
ordinary, by means of Neumann series. As in \cite{DeckKruse} we
define the two-dimensional convolution
\begin{equation} \label{conv}
f*g(t,x,s,y) \stackrel{def}{=} \int_s^t \ dr \int_R f(t,x,r,z) g(r,z,s,y) \ dz,
\end{equation}
then
\begin{equation} \label{Neuman}
 \Phi(t,x,s,y) =  \sum_{m=1}^{\infty} \Phi_m (t,x,s,y),
\end{equation}
where

$$
\Phi_2 (t,x,s,y) = \Phi_1(t,x,s,y) * \Phi_1 (t,x,s,y)
$$
and sequentially

$$
\Phi_{m+1} (t,x,s,y) = \Phi_1(t,x,s,y) * \Phi_m (t,x,s,y), \ m = 2,\ldots,
$$

so that the F.S. has the form
\begin{equation} \label{p representation}
p =Z + \sum_{m=1}^{\infty} Z * \Phi_m.
\end{equation}

\vspace{3mm}
Introduce also the function

\begin{equation} \label{Neuman2}
 \Phi^{(2)}(t,x,s,y) =  \sum_{m=2}^{\infty} \Phi_m (t,x,s,y),
\end{equation}

 \vspace{4mm}

 It remains to ground the convergence of the Neumann series \eqref{Neuman}. We will need some auxiliary estimates.

 \vspace{3mm}

 \ {\bf B.} \ Here and henceforth let $\delta$ be an arbitrary constant, \  $0 < \delta <
1$. The following estimate

$$
\left| \ \frac{\partial }{\partial x} \ Z_{t-s}(x-y) \ \right| \le (2 \pi e \delta)^{-1/2}\ (t-s)^{-1} \times
\exp \left\{ - \frac{(1 - \delta)(x-y)^2}{2(t-s)} \ \right\},
$$
holds, and consequently

\begin{equation} \label{Expr for Phi1 - 2}
\Phi_1(t,x,s,y) \le \frac{K[f]}{|x|^{\gamma}} \  (2 \pi e \delta)^{-1/2}\ (t-s)^{-1} \times
\exp \left\{ - \frac{(1 - \delta)(x-y)^2}{2(t-s)} \ \right\}.
\end{equation}

\vspace{3mm}

 \ {\bf C.} \ Let  $ \ A,B = \const > 0; \ $ then

\begin{equation} \label{Delta}
\min_{\Delta > 0} \left( \ A \ \Delta^{1 - \gamma} + B \Delta^{-\gamma} \ \right)  = \kappa(\gamma) \ A^{\gamma} \ B^{ 1 - \gamma},
\end{equation}
where

$$
\kappa(\gamma) =
\left(\frac{\gamma}{1-\gamma}\right)^{1-\gamma}+\left(\frac{\gamma}{1-\gamma}\right)^{-\gamma},
\ \ \gamma \in (0,1).
$$

 \vspace{3mm}

 \ {\bf D.} \ Recall the following identity which is well known in the probability theory, in the theory of Gaussian
 distribution:
\begin{equation} \label{Gauss conv}
\int_{\mathbb R} \exp \left\{ \  -0.5 \ \left(
\frac{(x-z)^2}{\sigma_1^2}  + \frac{(z-y)^2}{\sigma_2^2}  \ \right)
\ \right\} \ dz = (2\pi)^{1/2} \ \sigma_1 \ \sigma_2  \
\sigma_3^{-1} \ \exp \left\{ \ - 0.5 \frac{(x-y)^2}{\sigma_3^2} \
\right\},
\end{equation}
where $ \sigma_1, \ \sigma_2, \ \sigma_3  > 0$  satisfies
$$
\sigma_3^2 =  \sigma_1^2 + \sigma_2^2.
$$

\vspace{3mm}

 {\bf E.} \ Let us evaluate the following  important auxiliary parametric integral  which will be used after

\begin{equation} \label{int gamma}
J = J(x,t, y,s) := \int_{\mathbb R} |z|^{-\gamma} \ dz \ \exp
\left\{ \ -0.5 \left( \ \frac{(x-z)^2}{t-r}  +   \frac{(z-y)^2}{r-s}
\ \right) \ \right\}.
\end{equation}
 We split this integral $ J $ into two ones $J = J_1 + J_2$, where correspondingly

$$
  J_1 = J_1(\Delta) := \int_{|z| \le \Delta} (\cdot) \ dz, \ \ \ \  J_2 = J_2(\Delta): = \int_{|z| > \Delta} (\cdot) \ dz, \ \ \ \Delta = \const \in (0,\infty).
$$
 \ Evidently,

$$
J_1 \le 2 \int_0^{\Delta} z^{-\gamma} \ dz = \frac{2 \ \Delta^{1 - \gamma}}{1 - \gamma}.
$$
For the second integral we have
$$
J_2 \le \Delta^{-\gamma} \int_R  \exp \left\{ \ -0.5 \left( \ \frac{(x-z)^2}{t-r}  +   \frac{(z-y)^2}{r-s}  \ \right)     \ \right\} \ dz =: \Delta^{-\gamma} \ Y.
$$

Using equality \eqref{Gauss conv} we get

$$
Y =  (2 \pi)^{1/2} \ \sqrt{ (t-r)(r-s)(t-s)^{-1}} \ \exp \left\{  \ -0.5 \ \frac{(x-y)^2}{t-s} \    \right\},
$$
so that

$$
J \le \frac{2 \ \Delta^{1 - \gamma}}{1 - \gamma} + \Delta^{-\gamma} \ Y,
$$
and by \eqref{Delta}

\begin{equation}\label{aux int}
\begin{split}
J &\le \inf_{\Delta > 0} \left[ \frac{2 \ \Delta^{1 - \gamma}}{1 -
\gamma} + \Delta^{-\gamma} \ Y \right] = \kappa(\gamma) \  \left( \
\frac{2}{1-\gamma} \ \right)^{\gamma} \ Y^{1 - \gamma}\\
&=\kappa_1(\gamma)  \ \left[ (t-r)(r-s)/(t-s) \right]^{(1 -
\gamma)/2} \times \exp \left\{ \ - \frac{1-\gamma}{2} \
\frac{(x-y)^2}{t-s} \  \right\},
\end{split}
\end{equation}
where
$$
\kappa_1(\gamma) := (2 \pi)^{(1 - \gamma)/2}  \ \kappa(\gamma) \  \left( \ \frac{2}{1-\gamma} \ \right)^{\gamma}.
$$

\vspace{4mm}

 {\bf F.} \ Let us consider the following  more general parametric integral

\begin{equation} \label{int gamma  delta}
\ J_{\delta} = J_{\delta}(x,t, y,s) :=\int_{\mathbb R} |z|^{-\gamma}
\ dz \ \exp \left\{ \ -0.5(1 - \delta) \ \left( \
\frac{(x-z)^2}{t-r}  +   \frac{(z-y)^2}{r-s}  \ \right)     \
\right\},
\end{equation}
where $\delta = \const \in (0,1)$. We have, after a change of
variables,

$$
 J_{\delta}(x,t, y,s) = (1 - \delta)^{(\gamma - 1)/2} J(x\sqrt{1 - \delta},t, y \sqrt{1 - \delta},s),
$$
or equivalently

$$
J_{\delta}(x,t, y,s) = (1 - \delta)^{(\gamma - 1)/2} J(x,t/(1 - \delta), y,s/(1 - \delta)),
$$
and consequently

\begin{equation} \label{J delta}
J_{\delta}(x,t, y,s) \le (2 \pi)^{(1  - \gamma)/2} \ (1 -
\delta)^{(\gamma - 1)/2} \ \left[ (t-r)(r-s)/(t-s) \right]^{(1 -
\gamma)/2} \times\exp \left\{ \ - \frac{(1-\gamma)(1 - \delta)}{2} \
\frac{(x-y)^2}{t-s} \  \right\}.
\end{equation}

\vspace{4mm}

\ {\bf G.} The first term $ \ \Phi_1(\cdot) \ $ in  \eqref{Neuman}
is estimated in \eqref{Expr for Phi1 - 2}.
 Let us estimate the second one  $ \ \Phi_2. \ $   We have, after simple calculations,

$$
\Phi_2(x,t,y,s) \le \sqrt{1 - \delta} \ |x|^{-\gamma}
\frac{K^2[b]}{2 \pi e \ \delta} \ \int_s^t (t-r)^{-1} (r-s)^{-1} \
J(x \sqrt{1 - \delta}, t,y \sqrt{1 - \delta},s) \ dr,
$$

and, applying the estimate \eqref{aux int}, we have

$$
 \Phi_2(x,t,y,s) \le \kappa_1(\gamma)\, \sqrt{1 - \delta}  \ |x|^{-\gamma}  \times   \frac{K^2[b]}{2\pi e \ \delta} \
 \exp \left\{ - \frac{(1 - \gamma)(1 - \delta)}{2} \ \frac{(x-y)^2}{t-s}   \right\}
(t-s)^{(\gamma - 1)/2} \ I_2(t,s),
$$
where

\begin{equation*}
\begin{split}
I_2(t,s) & = \int_s^t (t-r)^{-(1+ \gamma)/2} \ (r-s)^{-(1+
\gamma)/2} \ dr = (t-s)^{-\gamma} \,B\left(\frac{1-\gamma}{2},\frac{1-\gamma}{2}\right)\\
&=(t-s)^{-\gamma} \,\frac{\Gamma^2((1 - \gamma)/2)}{\Gamma(1 -
\gamma)},
\end{split}
\end{equation*}
and $B(\cdot,\cdot)$ and $ \ \Gamma(\cdot) \ $ are respectively the
ordinary Beta and Gamma functions. Note that the last integral is
finite as long as $ \ \gamma \in [0,1). \ $
\par
 Evidently,

\begin{equation} \label{Phi2}
\Phi_2(x,t,y,s) \le \kappa_1(\gamma) \  |x|^{-\gamma} \ \sqrt{1 -
\delta} \  \frac{K^2[b]}{2\pi e \ \delta} \ (t-s)^{-(1+\gamma)/2} \
\, \frac{\Gamma^2((1 - \gamma)/2)}{\Gamma(1 - \gamma)} \
 \exp \left\{ - \frac{(1 - \gamma)(1 - \delta)}{2} \ \frac{(x-y)^2}{t-s}   \right\}.
\end{equation}

\vspace{4mm}

 \ {\bf H.} Let us investigate a more general case. Define the next functions

$$
G_1 = G_1(x,t,y,s) := (t-s)^{\alpha_1} \ |x|^{-\beta} \ \exp \left\{ \ -  0.5  \frac{(x-y)^2}{t-s}  \ \right\},
$$

$$
G_2 = G_2(x,t,y,s) := (t-s)^{\alpha_2} \ |x|^{-\gamma} \ \exp \left\{ \ -  0.5  \frac{(x-y)^2}{t-s}  \ \right\},
$$

where
$$
\alpha_1,\alpha_2 = {\const} > -1, \ \ \beta,\gamma = {\const} \in
[0,1), \ \ 0 < s < r < t < \infty,
$$

and define the convolution

$$
\ G = G(x,t,y,s) = G_1 * G_2(x,t,y,s) = \int_s^t dr \ \int_{\mathbb
R} G_1(x,t,z,r) \ G_2(z,r,y,s) \, dz .
$$

We  deduce analogously

\begin{equation} \label{G estim}
G(x,t,y,s) \le  \kappa_1(\gamma) \  (t-s)^{\frac{3}{2} + \alpha_1 +
\alpha_2 - \frac{\gamma}{2}} \  |x|^{-\beta}  \, B\left(\frac{3 -
\gamma}{2} + \alpha_1, \, \frac{3-\gamma}{2} + \alpha_2\right) \
\exp \left\{ \ - \frac{1-\gamma}{2} \ \frac{(x-y)^2}{t-s} \
\right\},
\end{equation}
so that
$$
B((3 - \gamma)/2 + \alpha_1, \ (3- \gamma)/2 + \alpha_2) = \frac{\Gamma((3 - \gamma)/2 + \alpha_1) \ \Gamma((3 - \gamma)/2 + \alpha_2)}{\Gamma(3 - \gamma + \alpha_1 + \alpha_2)}.
$$

\vspace{4mm}

 \ More generally, if  $ \ \delta = \const \in (0,1), \ $ then

$$
G_3^{(\delta)}(x,t,y,s) := \int_s^t dr \int_R dz \ |x|^{-\beta} \ |z|^{-\gamma} (t-r)^{\alpha_1} (r-s)^{\alpha_2} \times
$$

$$
\exp \left\{ - 0.5 (1 - \delta) \left[ \frac{(x-z)^2 }{t-r}  + \frac{(z-y)^2 }{(r-s)}  \right]  \right\} \le
$$

$$
  \kappa_1(\gamma) \ (1-\delta)^{-(1 + \beta)/2}   \  (t-s)^{3/2 + \alpha_1 + \alpha_2 + \gamma/2   } \  |x|^{-\beta} \  \times
$$

\begin{equation} \label{G delta estim}
B(3/2 + \alpha_1 + \gamma/2, 3/2 + \alpha_2 +\gamma/2) \ \exp \left\{ \  - \frac{(1 - \delta)(1-\gamma)}{2} \ \frac{(x-y)^2}{t-s} \  \right\}.
\end{equation}

\vspace{4mm}

 \ {\bf I.}  By means of induction and  using \eqref{G delta estim}
 we get

\begin{equation}
\begin{split}
\Phi_m(x,t,y,s) &\le (2 \pi)^{m/2 - 1} \ \kappa_1^{m-1}(\gamma) \
K^m[b] \ [(1 - \delta)(1 - \gamma)]^{m - 1} \ (e \
\delta)^{-m/2}\times
\\
&\times |x|^{-\gamma} \ (t-s)^{ 0.5m(1-\gamma) - 3/2 + \gamma/2  } \
\left[ \ \frac{\Gamma^m((1 - \gamma)/2)}{\Gamma(0.5 \ m \ (1 -
\gamma))} \ \right]\times\\
& \times \exp \left\{\ - 0.5 \ (1 - \delta) \ (1 - \gamma)^{m-1} \
\frac{(x-y)^2}{t-s} \ \right\},
\end{split}
\end{equation}

where the constant $ \ \gamma \in (0,1) \ $ is given, \ $  \ \delta
\in (0,1) \ $ is arbitrary. \par

\vspace{4mm}

 \ {\bf J.} Recall briefly the classical definition of the so-called (generalized) Mittag-Leffler function
 (\cite{Shukla}, \cite{Mittag})

$$
E_{\alpha,\beta}(z) \stackrel{def}{=} \sum_{k=0}^{\infty} \frac{z^k}{\Gamma(\alpha k + \beta)}.
$$

Here  $ \ \alpha = \const > 0, \ \beta \ge 0. \ $  We modify
slightly this definition. Namely, we introduce the function

\begin{equation} \label{Mitt}
g_{\mu}(z) \stackrel{def}{=} \sum_{m=1}^{\infty} \frac{z^m}{\Gamma(m
\ \mu)} = E_{\mu,0}(z) - 1, \ \  \mu = \const > 0.
\end{equation}

We need to use only real non-negative  values $z \in [0,\infty)$,
but the expression in \eqref{Mitt} is also defined for all the
complex values $ \ z. \ $
\par

The asymptotic as well as non-asymptotic behavior of these
functions, as $ \ |z| \to \infty \ $, has been investigated in
\cite{Shukla}, see also \cite{Evgrafov}. \par

We introduce also a slight modification of this function. Denote

\begin{equation} \label{psi}
\psi(z) = \psi_{\lambda,c,\mu}(z) \stackrel{def}{=} \sum_{m=1}^{\infty} \frac{c^m \ \exp \{ -\lambda^m \ z \}}{\Gamma(m \ \mu)}.
\end{equation}
and alike

\begin{equation} \label{psi2}
\psi^{(2)}(z) = \psi^{(2)}_{\lambda,c,\mu}(z) \stackrel{def}{=} \sum_{m=2}^{\infty} \frac{c^m \ \exp \{ -\lambda^m \ z \}}{\Gamma(m \ \mu)}.
\end{equation}

Here $ \ c,\lambda, \mu = \const > 0, \ $ but we need only the
values $ \ \lambda \ $ in the set $ \ (0,1). \ $ \par

 Evidently,
$$
\psi_{\lambda,c,\mu}(z) = \sum_{k=0}^{\infty} \frac{(-1)^k \
z^k}{k!} \cdot \left( \ E_{\mu,0}(c \ \lambda^k) - 1 \ \right)
=\sum_{k=0}^{\infty} \frac{(-1)^k \ z^k}{k!} \cdot g_{\mu} \left(\ c
\ \lambda^k \ \right).
$$


In order to obtain the asymptotic behavior of this function as $ \ z
\to \infty, \ $ we will use the methods in the book of V.N. Sachkov
(\cite[chapters 2,3]{Sachkov}), which is in turn the discrete analog
of the classical saddle-point method.
\par

Indeed, denoting $ \ \Lambda = |\ln(1 - \lambda)|, \ \lambda \in
(0,1) \ $, when $ z \ge e^e \Lambda \ $

\begin{equation} \label{Sachkov upp}
\psi_{\lambda,c,\mu}(z) \le S_1(\lambda,c,\mu) \cdot z^{C_1} \cdot
(\Lambda z)^{-(\mu/\Lambda) \cdot \, \ln (\ln (\Lambda z)/\Lambda)}
\end{equation}
 and alike

\begin{equation} \label{Sachkov low}
\psi_{\lambda,c,\mu}(z) \ge S_2(\lambda,c,\mu) \cdot z^{C_2} \cdot
(\Lambda z)^{-(\mu/\Lambda) \cdot \, \ln (\ln (\Lambda z)/\Lambda)
}.
\end{equation}

\vspace{4mm}

\begin{remark}
{\rm Notice that the function $\psi = \psi(z)$, introduced in
\eqref{psi}, decreases rapidly at $ z \to \infty$. Namely, let $B$
be an arbitrary positive constant, then
\begin{equation} \label{psi estim}
\exists C_1 = C_1(\lambda, c,\mu;B) \ : \ \psi_{\lambda,c,\mu}(z)
\le C_1(\lambda, c,\mu;B) (1 + z)^{-B}, \ \ \  \ z/\Lambda \ge e^e .
\end{equation}

\vspace{4mm}

 \ Further,  let us choose  $ \ \lambda= \lambda(\gamma) = 1 - \gamma, \ $

$$
 c = c(\delta,\gamma,K,\kappa_1) = \sqrt{2 \pi} \ \kappa_1 \ K \ [(1 - \delta) (1 - \gamma)] \ (e \ \delta)^{-1/2},
$$

$$
\mu = 1/(1- \gamma), \ z = z(\delta,\gamma,x,y,t,s) := 0.5 \ \frac{1- \delta}{(1 - \gamma)^{-1}} \ \frac{(x - y)^2}{t-s}
$$
and

\begin{eqnarray*}
\begin{split}
 D &= D(\kappa_1(\gamma),\delta,\gamma; x,y,t.s) := C_5(\kappa_1,\delta, \gamma) \ |x|^{-\gamma} \ (t-s)^{-(3 - \gamma)/2}
 \\
&=(2 \pi)^{-1} [\kappa_1(\gamma)]^{-1} \ [(1-\delta) \ (1 -
\gamma)]^{-1} \ |x|^{-\gamma} \ (t-s)^{-(3 - \gamma)/2},
\end{split}
\end{eqnarray*}
where evidently

$$
C_5 = C_5(\kappa_1,\delta, \gamma) =  (2 \pi)^{-1} [\kappa_1(\gamma)]^{-1} \ [(1-\delta) \ (1 - \gamma)]^{-1}.
$$

It follows immediately from this definition and from the relation
\eqref{Neuman} that

 $$
   \Phi(t, x, s,y) \le \overline{\Phi}(t,x,s,y), \
 $$
where

$$
\overline{\Phi}(t, x, s,y) \stackrel{def}{=} \inf_{\delta \in (0,1) } \overline{\Phi}_{\delta}(t, x, s,y),
$$

\begin{equation} \label{Event  Phi1}
\begin{split}
& \overline{\Phi}_{\delta}(t, x, s,y)  := \left[ \
D(\kappa_1(\gamma),\delta,\gamma; x,y,t.s)  \ \psi_{\lambda,c,\mu}
(z(\delta,\gamma,x,y,t,s)) \ \right]\\
&=(2 \pi)^{-1} [\kappa_1(\gamma)]^{-1} \ [(1-\delta) \ (1 -
\gamma)]^{-1} \ |x|^{-\gamma} \ (t-s)^{-(3 - \gamma)/2} \times
\psi_{\lambda,c,\mu} \left\{ 0.5 \ \frac{1- \delta}{(1 -
\gamma)^{-1}} \ \frac{(x - y)^2}{t-s} \ \right\}.
\end{split}
\end{equation}

}

\end{remark}

\vspace{4mm}

\ {\bf K.}  \ Ultimately, let now estimate the F.S. $ \ p = \
p(t,x,s,y) \ $ for the source equation \eqref{equation1}.
 Define again the following variables
$$
C_3 = C_3(\kappa_1, \delta,\gamma):= (2 \pi)^{-1} \ \kappa_1^{-1} \
(1 - \delta)^{-3/2} \ (1 - \gamma)^{-1} \ \Gamma((1 - \gamma)/2),
$$

$$
C_4  = C_4(\kappa_1, \delta,\gamma):= (2\pi)^{1/2} \ \kappa_1 \ [(1
- \delta)(1 - \gamma)] \ (e \delta)^{-1/2} \ \Gamma((1 - \gamma)/2),
$$
then, from the relations \eqref{Urav pt}, \eqref{Zfun}, \eqref{Event
Phi1}, we have
\begin{equation}\label{NasEstim1}
\begin{split}
& p(t,x,s,y) \le (2 \pi (t-s))^{-1/2} \ \exp \left\{ \  -
\frac{1}{2} \ \frac{(x-y)^2}{(t-s)}  \ \right\}\\
& + C_3 \times \sum_{m=1}^{\infty} C_4^m \  \frac{(t-s)^{0.5 m (1 -
\gamma)}}{ \Gamma(0.5 (m(1 - \gamma) + 2 - \gamma)   )  } \ \exp
\left\{ \ - 0.5 (1 - \delta)(1 - \gamma)^{m-1} \ \frac{(x -
y)^2}{t-s} \ \right\}.
\end{split}
\end{equation}

\vspace{4mm}
 The last expression may be simplified as follows. Suppose
$$
  (x-y)^2 \ge B \ \Lambda \ e^e \ (t-s), \ \ t > s,  \ \ B = \const > 0,
$$
the opposite case is trivial. Define a new function

\begin{equation} \label{qdelta fun}
\begin{split}
& q_{\delta}(t,x,s,y) := (2 \pi (t-s))^{-1/2} \ \exp \left\{ \  -
\frac{1}{2} \ \frac{(x-y)^2}{(t-s)}  \ \right\}\\
& +C_1(\lambda, c,\mu;B) \ C_5(\kappa_1,\delta, \gamma)  \
(t-s)^{-(3-\gamma)/2} \ \ \left[ \ 1 + \frac{1 - \delta}{(1 -
\gamma)^{-1}} \ \frac{(x-y)^2}{t-s} \  \right]^{-B},
\end{split}
\end{equation}
and

\begin{equation} \label{q fun}
q(t,x,s,y) := \inf_{\delta \in (0,1)}  q_{\delta}(t,x,s,y).
\end{equation}

 \vspace{4mm}

We deduce, under the above conditions, the following non-asymptotic
{\it upper estimate} for the fundamental solution of the equation
\eqref{equation1}

\begin{equation} \label{NasEstim2}
p(t,x,s,y) \le  q(t,x,s,y).
\end{equation}

\vspace{4mm}

Now we give the {\it lower estimate.} Suppose in addition that

$$
 K_- = K_-[b] := \inf_{t,x} b(t,x) > 0.
$$

As long as

$$
p \ge Z + c_1(K_-) Z * \Phi_1 +  c_2(K_-) Z *\Phi_2  - \sum_{m=3}^{\infty}Z*\Phi_m,
$$
we deduce, after some computations and under the same conditions,
that for {\it some} positive value  $ \ B_1 \ $

\begin{equation} \label{LowNasEstim3}
\begin{split}
& p(t,x,s,y) \ge (2 \pi (t-s))^{-1/2} \ \exp \left\{ \  -
\frac{1}{2}
\ \frac{(x-y)^2}{(t-s)}  \ \right\} \\
&+  C_6(\lambda, c,\mu;B_1) \ (t-s)^{-(3-\gamma)/2} \ \ \left[ \ 1 +
\frac{1 - \delta}{(1 - \gamma)^{-1}} \  \frac{(x-y)^2}{t-s} \
\right]^{-B_1}.
\end{split}
\end{equation}

\end{proof}

It is interesting to note that, for the parabolic PDE with smooth
coefficients, the bilateral estimates for the Fundamental Solution
has a Gaussian form (see, e.g.
 \cite{Aronson}, \cite{Eydelman2}).

\vspace{5mm}

\section{Solution of the equation \eqref{equation1} by means of the fundamental solution. }

\vspace{5mm}

\begin{theorem}
 Suppose, in addition to the conditions of Theorem \ref{Th FS}, that the continuous initial function
 $ f = f(x) $ in the
equation \eqref{equation1} has no more than polynomial growth at
infinity:

\begin{equation} \label{pol growth}
\exists \ W_0, \ W  \in (0, \infty)  \ : \ |f(x)| \le W_0 (1 +
|x|)^W, \  \ x \in \mathbb R.
\end{equation}

 Then the (unique) solution $ u = u(t,x) $  of the equation \eqref{equation1}, with initial Cauchy condition
 \eqref{initcond},
has the representation

\begin{equation} \label{representation}
u(t,x) = \int_{\mathbb R} p(t,x, s,y) \ f(y) \ dy,
\end{equation}
\end{theorem}
where $p$ is the F.S. of \eqref{equation1}.
 \vspace{4mm}

\begin{proof}
The relation \eqref{funsol} is obvious. The (uniform) convergence of
the integral in \eqref{p representation} follows immediately from
the condition \eqref{pol growth} and estimation \eqref{NasEstim1},
as well as the initial condition \eqref{initcond}. \par

Further, it is easy to verify, analogously to the proof of Theorem
\ref{Th FS}, that the convergence series in \eqref{p representation}
hods true also for the  partial derivatives $ \ \partial p/\partial
t,  \ \partial p/\partial x, $ and for  $ \ \partial^2 p/\partial
x^2. $\par

 See also \cite{Krylov1}, \cite{Lieberman}, \cite{Nash}. \par

This completes the proof.
 \end{proof}

\vspace{5mm}

\section{ Concluding remarks.}

\vspace{5mm}

 \ {\bf A.} The method here presented perhaps may be generalized on the equations of the form

\begin{equation} \label{equation slowly}
\frac{\partial u}{\partial t} = \frac{1}{2} \frac{\partial^2 u}{\partial x^2} + \frac{b(t,x) \ M(x)}{|x|^{\gamma}} \ \frac{\partial u}{\partial x},
\end{equation}
where $ \ M = M(x), \ x \in \mathbb R  \ $ is some continuous slowly
varying at the origin function.\par

\vspace{5mm}

 \ {\bf B.} The previous method may be extended to the more general $d$-dimensional case, $ d = 2,3,\ldots \ $, of the
parabolic differential equation of the form

$$
\frac{\partial u}{\partial t} = 0.5 \sum_{i=1}^d\ \sum_{j = 1}^d a_{i,j}(t,\vec{x}) \frac{\partial^2 u}{\partial x_i \partial x_j} +
\sum_{i=1}^d \frac{b_i(t,\vec{x})}{|x|^{\gamma(i)}} \ \frac{\partial u}{\partial x_i} + c(t,\vec{x}) u,
$$
$ \ u = u(t,\vec{x}),  \ $ with ordinary Cauchy initial condition

$$
\lim_{t \to s+} u(t,\vec{x})  = f(x), \ \ 0 \le s < t,  \ \ \vec{x}
= \{x_1, x_2,\ldots, x_d \} \in \mathbb R^d,
$$
where all the coefficients $ \  \{ a_{i,j}(\cdot,\cdot), \ b_i(\cdot,\cdot), \ c_i(\cdot,\cdot))  \}  \ $ are H\"older continuous,
the matrix $ \  \{ a_{i,j}(\cdot,\cdot) \} \ $ is symmetric, positive definite and bilateral bounded:

$$
\exists C_1, C_2 , \ 0 <  C_1 \le C_2 < \infty \ :  \ \forall
\vec{\xi} = \{ \xi_1,\xi_2, \ldots,\xi_d \} \ \in \mathbb R^d \
\Rightarrow
$$

$$
C_1 \sum_{i=1}^d \xi^2_i \le \sum_{i=1}^d \ \sum_{j = 1}^d a_{i,j}(t,\vec{x}) \xi_i \ \xi_j \le C_2\sum_{i=1}^d \xi^2_i;
$$

$$
\gamma(i) = {\const} \in (0,1),  \ \ f(\cdot) \in C(\mathbb R^d).
$$

\vspace{5mm}

 \ {\bf C.}  An open and interesting, by our opinion, problem: investigate the existence of F.S. and its properties for the initial problem for the
 parabolic differential equation with  variable coefficients of the (very singular) form

\begin{equation} \label{equation22}
\frac{\partial u}{\partial t} = \frac{1}{2} \frac{\partial^2 u}{\partial x^2} + \frac{b(t,x)}{|x|} \ \frac{\partial u}{\partial x} =: L_B[u].
\end{equation}

 The operator $ \ L_B \ $ is the infinitesimal operator for the so-called Bessel's non-homogeneous random process. \par

 The case $ \ b(t,x) = \const \ $ was considered in  \cite{Molchanov}, where was shown some probabilistic applications.
 A particular cases was represented in the monograph \cite{Revuz}, see also  \cite{Shiga}. \par
 \ A very interest application of the multivariate Bessel processes is described in the recent article \cite{Voit}, see also an article \cite{Diejen}. \par

\vspace{5mm}

 \ {\bf D.} Perhaps, one can study also the quasi-linear parabolic PDE of the form, for instance,

\begin{equation} \label{quasequation1}
\frac{\partial u}{\partial t} = \frac{1}{2} \frac{\partial^2 u}{\partial x^2} + \frac{b(t,x,u)}{|x|^{\gamma}} \ \frac{\partial u}{\partial x},
\end{equation}
with initial condition (\ref{initcond}),  or more generally

\begin{equation} \label{Qequat}
\frac{\partial u}{\partial t} = Q[t,\vec{x}, \ u, \ {\grad} \ u] \
u,
\end{equation}
where
\begin{equation} \label{def Q}
\begin{split}
 & Q[t,\vec{x}, \ u, \ {\grad} \ u] \ u \stackrel{def}{=} \ \frac{1}{2} \sum \sum_{i,j = 1}^d a_{i,j}(t,\vec{x}, u, \ {\grad} \ u) \ \frac{\partial^2 u}{\partial x_i \partial x_j}
 \\
 &+ \sum_{i=1}^d \frac{b_i(t, \vec{x},u, \ {\grad}\
u)}{|x_i|^{\gamma(i)}} \ \frac{\partial u}{\partial x_i} +
c(t,\vec{x}, u, \ {\grad} \ u)  u + h(t,\vec{x}),
\end{split}
\end{equation}

$$
 x = \vec{x} = \{ x_1,x_2,\ldots,x_d\}  \in \mathbb R^d, \ \ d = 1,2,3,\ldots, \  \ \gamma(i) \in [0,1),
$$
under appropriate conditions, e.g. with the symmetrical positive
definite  bounded and H\"older's continuous matrix $ \ \{
a_{i,j}(\cdot,\cdot,\cdot,\cdot)  \} \ $ as well as the coefficients
$ \ \{b_i(\cdot,\cdot,\cdot)\}, \ c(\cdot,\cdot,\cdot)  \  $ and
with the same initial condition \eqref{initcond}.\par

One can consider, for instance, following  the authors of  articles
\ \cite{Freidlin}, \ \cite{Inoue}, \ \cite{Korzeniowski}, \
\cite{Krylov2}, the  iteration procedure (recursion) consisting
only of the {\it linear equations}

$$
\frac{\partial u_{n+1}}{\partial t} = Q[t,\vec{x}, \ u_n, \ {\grad}
\ u_n] \, u_{n+1},  \ \ \ u_{n+1}(s+,x) = f(s,x), \ \ n =
0,1,2,\ldots
$$
with the initial condition relative to the number of iterations $n$

$$
\frac{\partial u_0}{\partial t} = Q[t,\vec{x}, \ 0, \  \vec{0}] \,
u_0, \ \ \  u_0(s+,x) = f(s,x).
$$

Evidently, if there exists the limit $u(t,x) := \lim_{n \to
\infty}u_n(t,x)$ in the space $ \ C_{\rm loc}^{1,2}(\mathbb
R_+,\mathbb R^d), \ $ it satisfies the equation \eqref{Qequat},
\eqref{def Q} with initial condition \eqref{initcond}.

\vspace{5mm}

 \vspace{0.5cm} \emph{Acknowledgement.} {\footnotesize The first
author has been partially supported by the Gruppo Nazionale per
l'Analisi Matematica, la Probabilit\`a e le loro Applicazioni
(GNAMPA) of the Istituto Nazionale di Alta Matematica (INdAM) and by
Universit\`a degli Studi di Napoli Parthenope through the project
\lq\lq sostegno alla Ricerca individuale\rq\rq }.\par

\end{document}